\documentclass{amsart}

\usepackage{amssymb}

\newtheorem{theorem}{Theorem}
\newtheorem{lemma}[theorem]{Lemma}
\newtheorem{proposition}[theorem]{Proposition}
\newtheorem{corollary}[theorem]{Corollary}
\newtheorem*{corthm}{Correspondence Theorem}

\theoremstyle{definition}

\DeclareMathOperator{\GL}{\mathsf{GL}}
\DeclareMathOperator{\supp}{supp}

\newcommand{\fa}[2]{\ensuremath {#1}\langle #2\rangle}
\newcommand{\fm}[1]{\ensuremath \langle #1 \rangle}
\newcommand{\matriz}[2]{\ensuremath M_{#1}(#2)}

\newcommand{\field}{\Bbbk}

\input tap.tex

\begin{document}

\title[Finitely generated invariants]{Finitely generated invariants of Hopf
  algebras on free associative algebras} 

\author[V. O. Ferreira]{Vitor O. Ferreira}
\address{Department of Mathematics - IME, University of S\~ao Paulo,
Caixa Postal 66281, S\~ao Paulo, SP, 05311-970, Brazil}
\email{vofer@ime.usp.br}
\thanks{
The first author was partially supported by CNPq, Brazil (Grant
302211/2004-7).}

\author[L. S. I. Murakami]{Lucia S. I. Murakami}
\address{Department of Mathematics - IME, University of S\~ao Paulo,
Caixa Postal 66281, S\~ao Paulo, SP, 05311-970, Brazil}%
\email{ikemoto@ime.usp.br}




\begin{abstract}
  We show that the invariants of a free associative algebra of
  finite rank under a linear action of a finite-dimensional Hopf algebra
  generated by group-like and skew-primitive elements form a finitely
  generated algebra exactly when the action is scalar. This generalizes
  an analogous result for group actions by automorphisms obtained by Dicks and
  Formanek, and Kharchenko.
\end{abstract}

\maketitle

\section*{Introduction}

Given a finite-dimensional vector space $V$ over a field $\field $ and a finite
subgroup $G$ of the group $\GL(V)$ of all invertible linear operators on $V$,
the action of the elements of $G$ on $V$ can be extended to the tensor algebra
$T(V)$ of $V$ in a natural way. The elements of $G$ become, then,
automorphisms of $T(V)$. We say that $G$ is a group of linear automorphisms of
the algebra $T(V)$. Fixing a basis, say, $\{x_1,\dots, x_r\}$, of $V$, $G$
can be regarded as a subgroup of the linear group $\GL(r,\field )$ and $T(V)$
becomes isomorphic to the free associative algebra $R=\fa{\field }{x_1,\dots, x_r}$.

The subalgebra $R^G$ of invariants of $R$ under the
action of $G$, defined to be the set of all elements $f\in R$ such that
$f^{\sigma} = f$, for all $\sigma\in G$, has been an object of interest for
some time. In particular, questions regarding presentations for it have been
addressed. Lane~\cite{dL76} and Kharchenko~\cite{vK78} have proved
independently that $R^G$ is a free algebra over $\field $ on a set of
homogeneous elements. Somewhat later, Dicks and Formanek~\cite{DF82}, and
Kharchenko~\cite{vK84}, again independently, showed that $R^G$ is a finitely
generated algebra exactly when $G$ is a group of scalar matrices, and is, therefore,
cyclic. Kharchenko's proof has been later simplified by Dicks and this new
argument appears in \cite[Theorem 10.4]{pC85} and in \cite[Theorem 32.7]{dP89}.

Lie algebras of derivations of free algebras have a similar
behavior with regards to constants. More precisely, given a finite-dimensional
restricted Lie algebra $L$ of linear derivations of $R$, 
Jooste~\cite{tJ78} and, independently, Kharchenko~\cite{vK81} have proved that
the subalgebra of constants $R^L = \{f\in R : f^{\delta} = 0 \text{, for
all $\delta\in L$}\}$ is free over $\field $ on a set of homogeneous
generators. It was then natural to consider the question regarding finite
generation of $R^L$. It follows from the work of Koryukin~\cite{aK94} that
exactly the same situation holds: $R^L$ is finitely generated as an algebra if
and only if $L$ consists entirely of scalar derivations. In \cite{FMta}, the
authors show that Dicks' proof for the automorphism case can be adapted to
take into account derivations. 

Actions of groups by automorphisms and of Lie algebras by derivations are
instances of Hopf algebra actions on rings. In fact, it was proved in
\cite{FMP04} that a free algebra under a homogeneous action by a Hopf algebra
has free invariants. In the present paper we address the question of finite
generation of the subalgebra of invariants. We show that whenever a
finite-dimensional Hopf algebra which is generated by group-like and
skew-primitive elements acts in a 
linear fashion on a free algebra, the subalgebra of invariants is a finitely
generated subalgebra if and only if the action is scalar. The proof of this
fact is based on Dicks' proof for the automorphisms case.

\section{Notation}

In this section we fix notation and nomenclature. Let $X$ be a non\-emp\-ty set and
let $\mathcal{F}=\fm{X}$ denote 
the free monoid on $X$. Let $\field $ be a field and let $R=\fa{\field }{X}$ denote the
free associative algebra on $X$ over $k$. The algebra $R$ is then a vector
space over $\field $ with basis $\mathcal{F}$. As usual, given
$f=\sum_{w\in\mathcal{F}} \lambda_w w\in R$,  with 
$\lambda_w\in \field $, for all $w\in\mathcal{F}$, the support of $f$ is defined to
be the following subset of $\mathcal{F}$:
$$
\supp(f) = \{w\in\mathcal{F} : \lambda_w\ne 0\}.\
$$

The free algebra $R$ can be graded by the usual degree function on $R$, that
is, $R = \bigoplus_{n\geq 0} R_n$, where $R_n$ stands for the linear span of
all the monomials of length $n$.

Let $H$ be a Hopf algebra and suppose that $H$ acts on $R$, that is, suppose
that $R$ is an $H$-module algebra (see \cite[Chapter 4]{sM93a}). We say that
the action of an element $h\in H$ on $R$ is 
\emph{linear} if it induces a linear
operator on the vector space $V=\sum_{x\in X}\field x$. In 
other words, the action of $h$ on $R$ is linear if for every 
$y\in X$, there exist scalars $\eta_{xy}(h)\in \field $, all but a finite number
of which nonzero, such that 
$$
h\cdot y = \sum_{x\in X} \eta_{xy}(h)x.
$$
Furthermore, we say that the action of $h$ on $R$ is \emph{scalar} if it is linear,
$\eta_{xx}(h) = \eta_{yy}(h)$, for all $x,y\in X$, and $\eta_{xy}(h) =
0$, for all $x,y\in X$ with $x\ne y$. The action of $H$ on $R$ is said to be
\emph{linear} if all $h\in H$ act linearly on $R$ and it is said to be
\emph{scalar} if it is linear and all the elements of $H$ act scalarly on $R$.
If $X$ is finite, say $X=\{x_1,\dots, x_r\}$, given $h\in H$, we often write
$\eta_{ij}(h)$ for $\eta_{x_ix_j}(h)$ and write $[h]_X$ for the matrix
$[\eta_{ij}(h)]\in\matriz{r}{\field }$. In this case, the action of an element $h\in
H$ on $R$ is scalar if it is linear and $[h]_X$ is a scalar matrix, say $[h]_X
= \eta I_r$,
for some $\eta\in \field $, where $I_r$ stands for the $r\times r$ identity
matrix. When this is the case we shall say that the action of $h$ is \emph{based on
$\eta$}.

If $H$ is a Hopf algebra with counit $\varepsilon$ and if $H$ acts linearly on
$R$, then the subalgebra of
invariants of $R$ under the action of $H$, defined by $R^H = \{f\in R : h\cdot
f = \varepsilon(h)f\text{, for all $h\in H$}\}$, is clearly a graded
subalgebra. Hence if $f\in R^H$, then all of the
homogeneous components of $f$ lie in $R^H$. 

Notation for Hopf algebra theory will be the usual (see \cite{mS69} or
\cite{sM93a}), including Sweedler's notation $\Delta(h) = \sum_{(h)}
h_{(1)}\otimes h_{(2)}$, for the comultiplication $\Delta$ on an element
$h\in H$.

Recall that an element $\sigma$ in a Hopf algebra $H$ is called a group-like
element if $\sigma\ne 0$ and  
$\Delta(\sigma)= \sigma\otimes\sigma$. The set of group-like elements is
linearly independent over $\field $ and forms a group under multiplication. 
Given group-like elements $\sigma,\tau\in H$, an element $\delta\in H$ is said to be 
$(\sigma,\tau)$-primitive if $\Delta(\delta) = \delta\otimes\sigma +
\tau\otimes \delta$. We say simply that $\delta$ is a skew-primitive element
if the reference to $\sigma$ and $\tau$ is not necessary. A primitive element
is just a $(1,1)$-primitive element. Finally if $\sigma$ is a group-like element, then
$\varepsilon(\sigma) =1$, whereas if $\delta$ is a skew-primitive element, then
$\varepsilon(\delta)=0$. (For these and other basic facts on Hopf algebras we
refer to \cite{mS69}.)

\section{Scalar actions}

We start by showing that finite generation of the invariants of a linear
action of a pointed Hopf algebra is a consequence of the action being scalar.

\begin{theorem} \label{th:fg}
  Let $\field $ be a field and let $H$ be a pointed Hopf $\field $-algebra which
  acts linearly on the free algebra $R = \fa{\field }{X}$ on a finite set $X$
  over $\field $. If the action of $H$ on $R$ is scalar, then $R^H$ is a finitely
  generated subalgebra. 
\end{theorem}

\begin{proof}
  By hypothesis, the action of each $h \in H$ is based on some $\lambda_h \in \field $. 
  Then, for $x,y \in X$,
  $$
  h \cdot(xy) = \sum_{(h)} (h_{(1)} \cdot x)(h_{(2)} \cdot y)  =
  \Big(\sum_{(h)} \lambda_{h_{(1)}} \lambda_{h_{(2)}} \Big) xy
  $$
  and, by induction on $n$, we have, for $x_1, x_2, \dots, x_n \in X$,
  $$
  h\cdot (x_1x_2 \dots x_n) =
  \Big(\sum_{(h)}\lambda_{h_{(1)}}\lambda_{h_{(2)}}\dots
  \lambda_{h_{(n)}}\Big)x_1x_2 \dots x_n.
  $$
  So, if $g \in \fa{\field }{X}$ is invariant and homogeneous of degree $t$,
  say $g = \sum_{i=1}^r \mu_iw_i$, where $\mu_i\in\field $ and $w_i$ are monomials
  of degree $t$, we have
  \begin{equation*}\begin{split}
  \varepsilon(h)g & = h \cdot g = \sum_{i=1}^r \mu_i\Big( \sum_{(h)} \lambda_{h_{(1)}}
  \dots \lambda_{h_{(t)}}\Big)  w_i \\
  & = \Big( \sum_{(h)} \lambda_{h_{(1)}}\dots
  \lambda_{h_{(t)}}\Big)g.
  \end{split}\end{equation*}
  Therefore, if there exists a homogeneous invariant element of degree $t$,
  then $\varepsilon(h) = \sum_{(h)} \lambda_{h_{(1)}}\dots 
  \lambda_{h_{(t)}}$, for all $h \in H$.
  It follows that every monomial of degree $t$ is invariant.
  
  If $R^H=\field $, then $R^H$ is trivially finitely generated. Otherwise, there
  exist invariant monomials of degree $\geq 1$. Let $t$ the least positive
  integer such that there exist an invariant monomial of 
  degree $t$. We shall show that $R^H$ is generated by the set of all monomials of degree $t$.
  Let $m>t$ and let $w = x_1\dots x_m$ be an invariant monomial of degree $m$.
  Write $m = qt+r$, with $0 \leq r < t$, and write $w = uv$, with $u$ and $v$
  monomials of degrees $qt$ and $r$, respectively.
  Observe that $u$ is invariant, because it is a product of $q$ monomials of
  degree $t$, all of which are invariant.
  Let $\{H_n\}$ be the coradical filtration of $H$. We shall proceed by
  induction on $n$.
  Since $H$ is pointed, $H_0 = \field G$, where $G$ stands for the set of all group-like
  elements of $H$. For each $\sigma \in G$, we have
  $$
  w = \varepsilon(\sigma)w = \sigma \cdot w = (\sigma \cdot
  u)(\sigma \cdot v) = u \lambda_{\sigma}^r v = \lambda_{\sigma}^r w.
  $$ 
  So $\lambda_{\sigma}^r = 1$. 
  Therefore, $\sigma \cdot v = \varepsilon(\sigma) v$, which means that $v$
  is invariant under the action of $H_0$. 
  Let $n>0$ and suppose that $v$ is invariant under the action of 
  $H_i$, for all $i = 0,1,\dots, n-1$. By the Taft-Wilson Theorem (see
  \cite[Theorem 5.4.1]{sM93a}), $H_n$ is generated by elements $h$ satisfying
  $\Delta(h) = h\otimes \sigma + \tau\otimes h + \sum l_i
  \otimes s_i$, with $\sigma,\tau \in G$ and $l_i,s_i \in H_{n-1}$. 
  Since $\varepsilon(h) = \sum_{(h)} \varepsilon(h_{(1)}h_{(2)})$, we have 
  $\varepsilon(h) = 2 \varepsilon(h) + \sum \varepsilon(l_is_i)$, so,
  $\sum \varepsilon(l_is_i) = -\varepsilon(h)$.
  Using the fact that $w, u\in R^H$ and $v$ is invariant under
  $H_i$, for $i<n$, we obtain
  \begin{equation*}\begin{split}
  \varepsilon(h)uv & = \varepsilon(h) w  = h \cdot w  \\
    & = (h \cdot u)(\sigma \cdot v) 
     + (\tau \cdot u)(h \cdot v) + \sum (l_i \cdot u)(s_i \cdot
     v) \\ 
    & = \varepsilon(h) w + u(h \cdot v) + \sum \varepsilon(l_i)
    \varepsilon(s_i) w \\ 
    & = \varepsilon(h) w + u(h \cdot v) - \varepsilon(h) w \\
    & = u(h\cdot v). 
  \end{split}\end{equation*}
  Hence $h \cdot v = \varepsilon(h)v$. We conclude that $v$ is an invariant
  monomial of degree $r$. By the minimality of $t$, it follows that $r=0$ and,
  thus, $m$ is a multiple of $t$.
 
  Therefore $R^H$ is generated by the set of all monomials of degree $t$ and, since
  $X$ is a finite set, $R^H$ is finitely generated.
\end{proof}

\section{Finitely generated invariants}

In the proof of \cite[Theorem 32.7]{dP89} and of the main theorem of
\cite{FMta}, essential use of the fact that there exist invariants and
constants, respectively, containing monomials with arbitrary initial segments
in their support is made. In each case, an appropriate power of a standard
polynomial is shown to possess such a property. This combinatorial ingredient
is no longer
at our disposal in the Hopf algebra context, for the action of the
symmetric group does not commute with the action of
skew-derivations. Nevertheless, we can resort to a substitute
construction resulting from the lemma below.  

We start by introducing some
notation. For each positive integer $n$, let
$c_n(Y,Z)$ be the following element of 
the algebra $\field [Y,Z]$ of commutative polynomials over the field $\field $ in the indeterminates
$Y$ and $Z$,
$$
c_n(Y,Z) = \sum_{i=0}^{n-1} Y^{n-1-i}Z^i.
$$
Note that $c_n(Y,Z)(Y-Z) = Y^n-Z^n$.

In what follows, $X$ will denote a nonempty set, $R=\fa{\field}{X}$ the free algebra on 
$X$ over $\field $, and $H$ a Hopf $\field $-algebra
which acts linearly on $R$. Given $w\in\mathcal{F}$, we shall write
$w\mathcal{F}$ for the set of
monomials in $\mathcal{F}$ that are left divisible by $w$, that is,
$$
w\mathcal{F} = \{wu : u\in\mathcal{F}\}.
$$

\begin{lemma}\label{le:jair}
  Suppose that $X$ is finite, say $X=\{x_1,\dots, x_r\}$, let
  $\sigma,\tau\in H$ be group-like elements and let $\delta\in H$ be a 
  $(\sigma,\tau)$-primitive element. Suppose that $\sigma$ and $\tau$ act scalarly on $R$,
  based on $\eta$ and $\mu$, respectively, and that $[\delta]_X$ is in Jordan normal
  form. Then for each positive integer
  $n$ and each $i=1,\dots, r$, there exists a nonzero $f\in R$ satisfying
  \begin{enumerate}
  \item $\delta\cdot f = c_n(\eta,\mu) g$, for some $g\in R$, and \label{c1le1}
  \item $\supp(f) \cap x_i\mathcal{F} \ne \emptyset$. \label{c2le1}
  \end{enumerate}
\end{lemma}

\begin{proof}
  Fix an $i=1,\dots, r$ and let $J$ be the Jordan block of $[\delta]_X$ containing the
  eigenvalue $\lambda$ on the $i$-th diagonal entry of $[\delta]_X$. Let $s$ be the
  positive integer such that the last occurrence of $\lambda$ in $J$ lies on
  the $s$-th diagonal entry of $[\delta]_X$. The figure below might help illustrate the
  choices of $J$ and $s$.

$$
\begintable
\begintableformat &\center \endtableformat
\B" | \- " " " " " " " " " " \- | \E" 
\B"^ | " " " " " " " " " " | \E"
\B" | " $\ddots$ " " " " " " " " " | \E" 
\B"_ | " " " " " " " " " " | \E"
\B" | " | @6\- | " " | \E" 
\B"^ | " | $\lambda$ " " " " " | " " | \E" 
\B" | " | $1$ " $\ddots$ " " " " | " " | \E"  
\B" $[\delta]_X = $ | " | " $\ddots$ " $\lambda$ " " " | " " | $\leftarrow \ i$ \E" 
\B" | " $J=$ | " " $1$ " $\ddots$ " " | " " | \E" 
\B"^ | " | " " " $\ddots$ " $\lambda$ " | " " | \E" 
\B"_ | " | " " " " $1$ " $\lambda$ | " " | $\leftarrow \ s$ \E" 
\B" | " | @6\- | " " | \E" 
\B"^ | " " " " " " " " " " | \E"
\B" | " " " " " " " " " $\ddots$ " | \E"
\B"_ | " " " " " " " " " " | \E" 
\B" | \- " " " " " " " " " " \- | \E" 
\B"^ " " " " " $\uparrow$ " " " $\uparrow$ " " " " \E" 
\B"^ " " " " " $i$ " " " $s$ " " " " \E" 
\endtable
$$

Then the element
  $$
  f = \sum_{j_1+\dots +j_n = i + (n-1)s} x_{j_1}\dots x_{j_n}
  $$
  satisfies both \eqref{c1le1} and \eqref{c2le1}. Indeed, it is easily checked
  that
  \begin{equation*}\begin{split}
  \delta\cdot (x_{j_1}x_{j_2}\dots x_{j_n})  = \ & \lambda c_n(\eta,\mu)
  x_{j_1}x_{j_2}\dots x_{j_n}  \\ 
  & \quad + \eta^{n-1}x_{j_1+1}x_{j_2}\dots x_{j_n}  \\
  & \quad + \eta^{n-2}\mu  x_{j_1}x_{j_2+1}\dots x_{j_n}\\
  & \quad + \dots \\
  & \quad + \eta\mu^{n-2}x_{j_1}\dots x_{j_{n-1}+1}x_{j_n} \\
  & \quad + \mu^{n-1} x_{j_1}\dots x_{j_{n-1}}x_{j_n+1},
  \end{split}\end{equation*}
  where each term in which there is an occurrence of $x_{s+1}$ should be
  interpreted as being zero. Therefore, $\delta\cdot f = \lambda c_n(\eta,\mu)
  f + f'$, where $f'$ is an element of $R$ with 
  \begin{equation*}\begin{split}
  \supp(f') \subseteq \{x_{k_1}\dots x_{k_n} : k_1+\dots&+ k_n = (i+1) + (n-1)s\\
  & \text{and } k_q \leq s \text{, for all $q=1,\dots, n$}
  \}.   
  \end{split}\end{equation*}
  Note that the restrictions $k_1+\dots+ k_n = (i+1) + (n-1)s$ and $k_q \leq
  s$, for all $q=1,\dots, n$, for the elements $x_{k_1}\dots x_{k_n}$ in
  $\supp(f')$ imply $k_q\geq 2$, for all $q=1,\dots, n$. Therefore, the
  element $x_{k_1}\dots x_{k_n}$ occurs in the supports of the image of the
  action of $\delta$ on $x_{k_1-1}x_{k_2}\dots x_{k_n}$, on $x_{k_1}x_{k_2-1}\dots
  x_{k_n}, \dots,$ and on $x_{k_1}x_{k_2}\dots x_{k_n-1}$ with coefficients
  $\eta^{n-1}, \eta^{n-2}\mu, \dots,$ and $\mu^{n-1}$,
  respectively\footnote{The authors thank Jair Donadelli Jr.~for this simple
  argument.}. 
  Hence, $f' =  c_n(\eta,\mu)f''$, for some $f''\in R$. It follows that
  $\delta\cdot f = c_n(\eta,\mu)(\lambda f+ f'')$. Thus, $f$ satisfies
  \eqref{c1le1}. For \eqref{c2le1}, note that,
  by definition of $f$, $x_ix_s\dots x_s\in\supp(f)$.
\end{proof}

For the next result, we need the following definition by Koryukin. Given non-negative
integers $i,j,k$, let $\tau_{ijk}$ be the
linear operator of $R_{i+j+k}$ which, on
monomials $u\in R_i$, $v\in R_j$, and $w \in R_k$, satisfies 
$\tau_{ijk}(uvw) = uwv$. A graded subalgebra $S$ of $R$ is called an
\emph{algebra with inserts} if $\tau_{ijk}(S_{i+j}S_k) \subseteq S_{i+j+k}$,
for all non-negative $i,j,k$. The next result is due to Koryukin. We include a
proof for the reader's convenience.  

\begin{proposition}[{\cite[Lemma~1.6]{aK94}}]
  The subalgebra of invariants $R^H$ is an algebra with inserts.
\end{proposition}

\begin{proof}
  Let $i,j,k$ be non-negative integers,
  let $u$ and $v$ be monomials of degrees $i$ and $j$,
  respectively, and let $g \in R^H_k$. Then for all $h\in H$, we have
  \begin{equation*}\begin{split}  
  h \cdot \tau_{ijk}(uvg) & = h \cdot(ugv) = \sum_{(h)}
  (h_{(1)}\cdot u)(h_{(2)} \cdot g)(h_{(3)} \cdot v) \\
  &  = \sum_{(h)} (h_{(1)} \cdot u) \varepsilon(h_{(2)}) g (h_{(3)}
  \cdot v) \\ 
  & =  \sum_{(h)}(h_{(1)} \cdot u) g (h_{(2)} \cdot v) =
  \tau_{ijk}((h\cdot(uv))g),
  \end{split}\end{equation*}
  since the action is linear.
  Hence, by the linearity of $\tau_{ijk}$, we have, for all $f \in
  R_{i+j}$ and $g \in R_k^H$,
  $$h \cdot \tau_{ijk}(fg) = \tau_{ijk}((h \cdot f)g).$$
  Therefore, if $f \in R^H_{i+j}$ and $g \in R^H_k$ then $\tau_{ijk}(fg)
  \in R^H$.
\end{proof}

\begin{corollary} \label{co:xl}
  Let $x \in X$. If there exists $f \in R^H$ with $\supp(f) \cap x
  \mathcal{F} \neq \emptyset$ then, for each positive integer $k$, there exists
  $\bar{f} \in R^H$ such that $\supp(\bar{f}) \cap x^k \mathcal{F}
  \neq \emptyset$. 
\end{corollary}

\begin{proof}
  We prove this fact by induction on $k$. We can assume that $f$ is
  homogeneous, for the action of $H$ is linear. Let $d$ be the degree of $f$
  and let $xw \in \supp(f)$. For $k=1$, there is nothing to prove. Suppose
  $k>1$ and assume that there exists a homogeneous $\tilde{f} \in R^H$ of degree $t$, with
  $\supp(\tilde{f}) \cap x^{k-1}\mathcal{F} \neq \emptyset$, say
  $x^{k-1}\tilde{w} \in \supp(\tilde{f})$.
  Then $\bar{f} = \tau_{k-1,t-k+1,d}(\tilde{f}f)$ is invariant, by the
  previous lemma, and $x^kw\tilde{w} \in \supp(\bar{f})$.
\end{proof}

The following is a well known fact. It
follows, for instance, from \cite[Lemma~5.5.1]{sM93a}.   

\begin{proposition}\label{prop:pointed}
  A Hopf algebra which is generated by group-like and skew-primitive
  elements is pointed. \qed
\end{proposition}


Finally, we shall make use of the following Galois correspondence.

\begin{corthm}[{\cite[Theorem~1.2]{FMP04}}]
  Let $\field $ be a field, let $X$ be a set with $|X|>1$, and let $R=\fa{\field}{X}$ 
  be the free algebra on $X$ over $\field $. Let $H$ be a
  finite-dimensional pointed Hopf $\field $-algebra which acts faithfully and linearly
  on $R$. Then there exists an inclusion-inverting one-to-one correspondence
  between the set of all free subalgebras of $R$ containing $R^H$ and the set
  of all right coideal subalgebras of $H$. \qed
\end{corthm}

We are ready to state and prove the main result of the paper.

\begin{theorem}\label{th:main}
  Let $\field $ be a field, let $X$ be a nonempty set, and let $R=\fa{\field }{X}$ be the free
  algebra on $X$ over $\field $. Let $H$ be a finite-dimensional Hopf $\field $-algebra
  which acts faithfully and linearly on $R$. Suppose that $H$ is generated by
  group-like and skew-primitive elements. Then the subalgebra of invariants $R^H$
  is finitely generated if and only if $X$ is finite and the action of $H$ on
  $R$ is scalar. 
\end{theorem}

\begin{proof}
  If $X$ is finite and the action of $H$ on $R$ is scalar then $R^H$
  is finitely generated by Theorem \ref{th:fg}. Conversely, suppose
  $R^H$  is finitely generated. Then there exist $x_1, \dots, x_n \in
  X$ such that $R^H \subseteq \fa{\field }{x_1,\dots,x_n}$.
  If $X$ were infinite there would exist an infinite chain 
  \begin{equation*}\begin{split}
  R^H & \subseteq \fa{\field }{x_1,\dots,x_n} \varsubsetneq
  \fa{\field }{x_1,\dots,x_n,x_{n+1}} \\
  & \varsubsetneq \dots \varsubsetneq
  \fa{\field }{x_1,\dots,x_n,\dots,x_{n+m}} \varsubsetneq \dots
  \end{split}\end{equation*}
  of free subalgebras of $R$ containing $R^H$. 
  By the Correspondence Theorem, there would also exist an infinite
  chain of subalgebras of $H$, which is impossible, since $H$ is
  finite-dimensional. Therefore, $X$ is finite.

  If $X$ contains one single element, then any linear action is scalar. Thus
  we can assume that $|X|>1$, say $X = \{x_1, \dots, x_r\}$.
  By Proposition~\ref{prop:pointed}, $H$ is pointed, so its coradical is $H_0 =
  \field G$, where $G$ is the set of all group-like elements of $H$. 
  By the Correspondence Theorem, $R^{H_0}$ is a free subalgebra
  of $R$ containing $R^H$ and, by the same argument as above,
  $R^{H_0}$ is finitely generated. Since $H$ is finite-dimensional, $G$ is
  finite. So, by \cite[Theorem~5.3]{DF82}, the action of
  $H_0$ on $R$ is scalar.
  
  Since $H$ is generated by $G$ and skew-primitive elements, it remains to
  show that each skew-primitive element acts scalarly on $R$. We start by
  remarking that it can be assumed that $\field$ is algebraically closed, for if
  $\bar{\field}$ 
  denotes the algebraic closure of $\field$, then the action of $H$ on $R$
  induces an action of $\bar{H} = \bar{\field}\otimes_{\field} H$ on
  $\bar{R} = \bar{\field}\otimes_{\field}R$ such that
  $\bar{R}^{\bar{H}} = \bar{\field}\otimes_{\field} R^H$ is
  finitely generated. 

  Let $\delta \in H$ be a $(\sigma,\tau)$-primitive element, where
  $\sigma,\tau$  are group-like elements with actions based on $\eta$ and
  $\mu$, respectively. Consider the subalgebra
  $H(\delta)$ of $H$ generated by $\{\delta,\sigma, \tau\}$. Then $H(\delta)$ is clearly a
  Hopf subalgebra of $H$ and, by the Correspondence Theorem, the subalgebra of
  invariants $R^{H(\delta)}$ of $R$ under $H(\delta)$ is finitely generated. Now, under
  our assumption that $\field$ is algebraically closed, we can further assume
  that $[\delta]_X$ is in Jordan
  normal form, since any basis for the vector space $V = \sum_{x\in X}\field x$
  is a set of free generators of $R$ giving rise to the same grading.

  We shall start by showing that $[\delta]_X$ is a diagonal matrix. Suppose otherwise. After
  reordering the basis, if necessary, we would have 
  $\delta \cdot x_1 = \lambda x_1 + x_2$ and 
  $\delta \cdot x_2 = \lambda x_2 + \zeta x_3$,
  where $\lambda,\zeta \in \field$, $\zeta = 0$ or $\zeta = 1$, and
  $x_1,x_2 \not\in \supp(\delta\cdot x_i)$, for $i \geq 3$.
  Let $A$ be the set of all monomials different from $1$ which occur in
  the support of the elements of a finite set of generators for $R^{H(\delta)}$.
  We claim that there is an $m \geq 1$ such that $x_1^m \in A$. In fact,
  it is sufficient to show that there exists $\phi \in R^{H(\delta)}$ satisfying
  \begin{equation}\label{inv}\tag{$\ast$}
  \supp(\phi)\cap x_1 \mathcal{F} \neq \emptyset
  \end{equation}
  Indeed, having such an invariant $\phi$, using Corollary~\ref{co:xl},
  we can produce $\bar{\phi} \in R^{H(\delta)}$ with $\supp(\bar{\phi}) \cap 
  x_1^k \mathcal{F} \neq \emptyset$, say $x_1^kw \in \supp(\bar{\phi})$,
  where $k$ is an integer greater than the degrees of the elements of
  $A$. Since $x_1^kw$ is a product of elements of $A$ we must have $x_1^m
  \in A$, for some $m \geq 1$.
 
  In order to exhibit such an invariant element, let $f \in R$ be the
  element obtained in Lemma \ref{le:jair} with
  $i=1$ and $n = |G|$. We have
  $\delta \cdot f = c_n(\eta,\mu)g$, for some $g\in R$, and $\supp(f) \cap
  x_1 \mathcal{F} \neq \emptyset$. Since $\sigma,\tau\in G$, we have
  $\sigma^n=\tau^n=1$ and, therefore, $\eta^n=\mu^n=1$. If, on the one
  hand, $\eta \neq \mu$, it follows that $c_n(\eta,\mu)=0$; hence $\delta \cdot f = 0 =
  \varepsilon(h)f$, that is, $f$ is invariant under $\delta$. Since $f$ is
  homogeneous of degree $n$ it also follows that $\sigma\cdot f = \tau\cdot f
  = f$. Thus $f\in R^{H(\delta)}$. So $\phi=f$ satisfies \eqref{inv}. On the
  other hand, if $\eta = \mu$, then $\sigma = \tau$, for the action of $H$ on 
  $R$ is faithful. In this case $\sigma^{-1}\delta$ is a primitive
  element. Therefore, the characteristic of $\field$ must be positive,
  otherwise, the subalgebra of $H$ generated by $\sigma^{-1}\delta$ would be
  infinite-dimensional. Since $c_n(\eta,\eta) = n\eta^{n-1}$, we have
  $\delta \cdot f^k = kn\eta^{n+k-2}g^k$, for every $k \geq 1$. If $p$ is the
  characteristic of $\field$, then $\delta
  \cdot f^p=0=\varepsilon(h)f^p$, that is, $f^p$ is invariant under
  $\delta$. Again, because it is homogeneous of degree $pn$, $f^p$ is also
  invariant under $\sigma$ and $\tau$. So, \textit{a fortiori,} $f^p\in R^{H(\delta)}$. In
  this case $\phi=f^p$ is an invariant element satisfying \eqref{inv}. 

  Now take a homogeneous $z \in R^{H(\delta)}$ of degree $m$ with $x_1^m \in
  \supp(z)$ and write
  $z = x_1^m + \nu x_1^{m-1}x_2+\bar{z}$, with $\nu \in \field$ and
  $x_1^m,x_1^{m-1}x_2 \not\in \supp(\bar{z})$.
  Then
  \begin{equation*} \begin{split}
   \delta \cdot z & = \lambda c_m(\eta,\mu)x_1^m +
   \mu^{m-1}x_1^{m-1}x_2 \\
   & \quad + \nu \lambda \eta c_{m-1}(\eta,\mu)
   x_1^{m-1} x_2 + \nu \lambda \mu^{m-1} x_1^{m-1}x_2 + \bar{\bar{z}}
  \end{split} \end{equation*}
  where $x_1^m, x_1^{m-1}x_2 \not\in \supp(\bar{\bar{z}})$ 
  Since $\delta \cdot z=0$, by comparing coefficients, we get
  $\lambda c_m(\eta,\mu)=0$ and 
  $0 = \nu\lambda(\eta c_{m-1}(\eta,\mu) + \mu^{m-1})
  + \mu^{m-1} = \nu \lambda c_m(\eta,\mu) + \mu^{m-1} =
  \mu^{m-1},$  which implies $\mu=0$, a contradiction.
  Therefore $[\delta]_X$  must be a diagonal matrix.
  So, there exist $\lambda_1, \dots, \lambda_r \in \field$ such that 
  $\delta\cdot x_i = \lambda_ix_i$, for $i = 1, \dots, r$.
  
  In order to show that all these $\lambda_i$ are equal, we observe
  that there exists $w \in \mathcal{F}$ such that $wx_i \in A$ for
  every $i = 1, \dots, r$, otherwise we would be able to construct a
  monomial $\bar{w}$ of degree $k$ (where $k$ was chosen to be an
  integer greater than the degrees of all the elements of $A$) whose initial
  segments would all 
  lie outside $A$. However this could not happen, since, for each
  $i=1,\dots, r$, there exists $f_i \in R^{H(\delta)}$ with $x_i \mathcal{F} \cap
  \supp(f_i) \neq \emptyset$ (this can be proved in the same way as we
  have done for $i=1$ using Lemma \ref{le:jair}). But because $R^{H(\delta)}$ is an
  algebra with inserts we would eventually be able to produce an $f \in R^{H(\delta)}$
  with $\bar{w} \mathcal{F} \cap \supp(f) \neq \emptyset$. Now this would
  imply that some initial segment of $\bar{w}$ should be an element of
  $A$; a contradiction.
  Write $w = x_{i_1}x_{i_2} \dots x_{i_t}$. So
  $\delta \cdot w = \xi w$, where $\xi = \sum_{j=1}^t
  \lambda_{i_j}\eta^{t-j}\mu^{j-1}$. 
  Therefore
  $0 = \delta \cdot (wx_i) = \xi \eta wx_i + \mu^t \lambda_iwx_i$
  so, $\lambda_i = -\mu^{-t}\eta\xi$, for every $i = 1,\dots,r$.
  Hence, the action of $\delta$ is scalar.
\end{proof}

We believe the hypothesis on $H$ being generated by group-like and
skew-primitive elements not to be too restrictive, for a great number of
finite-dimensional pointed Hopf algebras do have this property. In fact, it has
been conjecture in \cite[1.4]{AS00} that all finite-dimensional pointed Hopf
algebras over an algebraically closed field of characteristic $0$ are generated
as algebras by group-like and skew-primitive elements. In view of this fact it
would not seem unreasonable to believe that Theorem~\ref{th:main} holds for
arbitrary finite-dimensional pointed Hopf algebras, but we do not have a proof
for this fact.


\begin{thebibliography}{FMP04}

\bibitem{AS00}
  \textsc{N. Andruskiewitsch and H.-J. Schneider},
  Finite quantum groups and Cartan matrices,
  \textit{Adv. Math.}~\textbf{154} (2000), no. 1, 1--45.

\bibitem{pC85}
  \textsc{P. M. Cohn},
  \textit{Free Rings and Their Relations}, 2nd. Ed.,
  Academic Press, London, 1985.

\bibitem{DF82}
  \textsc{W. Dicks and E. Formanek},
  Poincar\'e series and a problem of S.~Montgomery,
  \textit{Linear and Multilinear Algebra}~\textbf{12} (1982/83), no. 1, 21--30.

\bibitem{FMta}
  \textsc{V. O. Ferreira and L. S. I. Murakami},
  Finitely generated constants of free algebras,
  to appear in: A. Giambruno, C. Polcino Milies, S.K. Sehgal (Eds.), Groups,
  Rings, and Group Rings, \textit{Lecture Notes in Pure and Appl. Math.},
  \textbf{247}, Dekker,  New York. 

\bibitem{FMP04}
  \textsc{V. O. Ferreira, L. S. I. Murakami and A. Paques},
  A Hopf-Galois correspondence for free algebras,
  \textit{J. Algebra}~\textbf{276} (2004), no. 1, 407--416.

\bibitem{tJ78}
  \textsc{T. W. Jooste},
  Primitive derivations in free associative algebras,
  \textit{Math. Z.}~\textbf{164} (1978), no. 1, 15--23.

\bibitem{vK78}
  \textsc{V. K. Kharchenko},
  Algebras of invariants of free algebras,
  \textit{Algebra i Logika}~\textbf{17} (1978), no. 4, 478--487, 491.
  (English translation in \textit{Algebra and Logic}~\textbf{17} (1978),
  no. 4, 316--321 (1979).) 

\bibitem{vK81}
  \bysame,
  Constants of derivations of prime rings,
  \textit{Izv. Akad. Nauk SSSR Ser. Mat.}~\textbf{45} (1981), no. 2, 435--461,
  464. (English translation in 
  \textit{Math. USSR-Izv.}~\textbf{18} (1982), no. 2, 381--401.)

\bibitem{vK84}
  \bysame,
  Noncommutative invariants of finite groups and Noetherian varieties,
  \textit{J. Pure Appl. Algebra}~\textbf{31} (1984), no. 1-3, 83--90.

\bibitem{aK94}
  A. N. Koryukin,
  On noncommutative invariants of bialgebras,
  \textit{Algebra i Logika}~\textbf{33} (1994), no. 6, 654--680, 716. (English
  translation in \textit{Algebra and Logic}~\textbf{33} (1994), no. 6, 366--380
  (1995).)

\bibitem{dL76}
  \textsc{D. Lane},
  \textit{Free Algebras of Rank Two and Their Automorphisms},
  PhD thesis, London, 1976.

\bibitem{sM93a}
  \textsc{S. Montgomery},
  \textit{Hopf Algebras and Their Actions on Rings},
  Amer. Math. Soc., Providence, RI, 1993.

  
\bibitem{dP89}
  \textsc{D. S. Passman},
  \textit{Infinite Crossed Products},
  Academic Press, Boston, MA, 1989.

\bibitem{mS69}
  \textsc{M. Sweedler},
  \textit{Hopf Algebras},
  W. A. Benjamin, New York, 1969.

\end{thebibliography}
\end{document}